\documentclass[11pt]{article}
\usepackage{amssymb,amsmath,amsfonts,latexsym} 
\newtheorem{lemma}{Lemma}[section]
\newtheorem{thm}[lemma]{Theorem}

\newtheorem{corollary}[lemma]{Corollary}
\newtheorem{proposition}[lemma]{Proposition}
\newtheorem{fact}[lemma]{Fact}

\newtheorem{conjecture}[lemma]{Conjecture}

\newcommand{\comment}[1]{}

\def\binom#1#2{{#1\choose#2}}
\def\ex{{\text{\rm ex}}}

\def\G33{{{\mathbb G}_{3\times 3}}}
\def\I2{{\mathbb I}_{\geq 2}}
\def\M5{{\mathbb M}_5}
\def\P4{{\mathbb P}_4}
\def\qed{\ifhmode\unskip\nobreak\hfill$\Box$\medskip\fi\ifmmode\eqno{\Box}\fi}

\def\cE{{\mathcal E}}
\def\bbF{{\mathbb F}}
\def\bF{{\mathbf  F}}
\def\cF{{\mathcal F}}

\def\bbG{{\mathbb G}}

\def\bbH{{\mathbb H}}

\def\cH{{\mathcal H}}

\def\cP{{\mathcal P}}

\def\bbT{{\mathbb T}}

\begin{document}

\pagestyle{myheadings}
\markright{{\small \sc Z. F\"uredi:}
  {\it\small $2$-cancellative families}}
\thispagestyle{empty}

\title{\Huge $2$-cancellative hypergraphs and codes
\footnote{ This copy was printed on {\today}.\quad
    {\rm\small {\jobname}.tex,} \hfill Version as of 
    March 9, 2011.
\break\indent{\it Keywords:} Tur\'an hypergraph problems, 
cancellative hypergraphs, superimposed codes, independent polynomials.
\hfill \break\indent{\it 2010 Mathematics Subject Classifications:} 05D05, 11T06, 05D40.} }
\author{{\bf Zolt\'an F\"uredi}
\thanks{ Research supported in part by the Hungarian National Science Foundation
 OTKA, and by the National Science Foundation under grant NFS DMS 09-01276.}
\\ Department of Mathematics, University of Illinois at Urbana-Champaign,
\\ Urbana, IL 61801, USA \quad and
\\ R\'enyi Institute of Mathematics of the Hungarian Academy of Sciences,
\\ Budapest, P. O. Box 127, Hungary-1364
\\ e-mail: {\tt z-furedi@illinois.edu}\quad  and\quad {\tt furedi@renyi.hu}
}

\date{${}$}
\maketitle
\begin{abstract}
A family of sets $\mathcal F$ (and the corresponding family of 
0-1 vectors) is called $t$-{\bf cancellative} 
if for all distict $t+2$ members  $A_1, \dots, A_t$ and $B,C\in \cF$ 
  $$A_1\cup\dots \cup A_t\cup  B \neq A_1\cup \dots A_t \cup C. 
  $$
Let $c(n,t)$ be the size of the largest $t$-cancellative family on $n$
elements, and let $c_k(n,t)$ denote the largest $k$-uniform family.
We significantly improve the previous upper bounds, e.g., we show
 $c(n,2)\leq 2^{0.322n}$ (for $n> n_0$). 
Using an algebraic construction we show that 
 $c_{2k}(n,2)=\Theta (n^k)$ for each $k$ when $n\to \infty$. 

\end{abstract}

\section{Introduction, definitions}

There are many instances in Coding Theory when codewords must 
be restored from partial information like defected data (error 
correcting codes) or some superposition of the strings (these 
can lead to Sidon sets, sum-free sets, etc). 
A family of sets $\mathcal F$ (and the corresponding family of 
0-1 vectors) is called {\bf cancellative} if $A$ and $A\cup B$ 
determine $B$ (in case of $A, B\in \mathcal F$ and $A\neq A\cup B$). 
For a precise definition we require that for all $A,B,C\in \mathcal F$, $A\neq B$, $A\neq C$
  $$A\cup B=A\cup C\enskip  \Longrightarrow\enskip  B=C. $$
Let $c(n)$ be the size of the largest cancellative family on $n$
elements, and let $c_k(n)$ denote the size of the largest $k$-uniform family.
This definition can be extended (see above in the abstract).
In this paper we focus on $2$-cancellative $r$-uniform hypergraphs,
 i.e., families of $r$-sets, and  on $2$-cancellative {\it codes}, 
 where there is no restriction on the sizes of the hyperedges. 

Speaking about a  hypergraph $\mathbb F=(V, {\mathcal F})$ we frequently
 identify the vertex set $V=V(\bbF)$  by the set of first integers $[n]:=\{ 1,2, \dots, n\}$, 
 or elements of a $q$-element finite field $\bF_q$.
To shorten notations we frequently say  'hypergraph $\cF$' (or set system $\cF$) thus
 identifying $\bbF$ to its edge set $\cF$. 
The {\em degree}, $\deg_\bbF(x)$,  of an element $x\in [n]$ is the number of hyperedges 
 of ${\cal F}$ containing $x$.
The hypergraph $\cF$ is {\em uniform} if every edge has the same number of elements,
 $r$-uniform means $|F|=r$ for all $F\in \cF$.
An $r$-uniform hypergraph $(V, \cF)$ is called $r$-{\em partite} if 
 there exists an $r$-partition of $V$, $V=V_1\cup \dots \cup V_r$, such that 
 $|F\cap V_i|=1$ for all $F\in \cF$, $i\in [r]$.

Let $f(n, P_1, P_2, \dots )$ denote the  maximum number of
 subsets which can be selected from $\{1,\dots ,n\}$ satisfying all the properties
 $P_1, P_2, \dots$.
With this notation $c(n,t):= f(n,t\text{\rm{-CANC}})$, where $t$-CANC stands for 
 $t$-cancellativeness.

A hypergraph is {\em linear}, if $|E\cap F|\leq 1$ holds for every pair of edges.
An $(n,r,2)$-{\em packing}  is a linear $r$-uniform hypergraph $\cP$ on $n$ vertices.
Obviously, $|\cP| \leq   \binom{n}{2} / \binom{r}{2}$.
If here equality holds, then $\cP$ is called an $S(n,r,2)$  {\em Steiner system}.

\section{Cancellative and locally thin families}

The asymptotics of the maximum size of a cancellative family was given by  
  Tolhuizen~\cite{T} (construction) and in~\cite{Fur27} (upper bound) that
there exists a $C>0$ such that 
\begin{equation*}  
    \frac{C}{\sqrt{n}} \, 1.5^n< c(n) < 1.5 ^n
   \end{equation*}
  holds. 
The problem was proposed by Erd\H os and Katona~\cite{Kat} who conjectured 
 $c(n)=O(3^{n/3})$ which was disproved by an elegant construction
  by Shearer showing $c(3k)\geq k3^{k-2}$ leading to $c(n)> 1.46^n$ for $n> n_0$. 
Since a kind of product of two cancellative families is again cancellative
 we have $c(n+m)\geq c(n)c(m)$.
Thus the $\lim c(n)^{1/n}$ obviously exists.
This is not known for $2$-cancellative hypergraphs, so
 K\"orner and Sinaimeri~\cite{KS} introduced
  $$ t(4):=\limsup_{n\to \infty}\frac{1}{n}\log_2 c(n,2) $$ and proved
  $0.11<t(4)\leq0.42$. 
As usual all logarithms have base two. 
The lower bound follows from a standard probabilistic argument.
We will show that $t(4)\leq \log_2 5 -2= 0.3219\dots$

\begin{thm}\label{th:2c} \quad
$c(n,2)< 9\sqrt{n} \left( \dfrac{5}{4}\right)^n.$ 
 \end{thm}

The proof is postponed to Section~\ref{ss:2c}. 
Without loss of generality we can suppose that the $n$-set underlying $\cF$ is $[n]$.
We associate to every subset $A\in \cF$, its characteristic binary vector, 
 ${\mathbf x}:={\mathbf x}(A) = (x_1,\dots, x_n)$,
 with $x_i = 1$ if $i\in A$ and $x_i = 0$ otherwise. 
One can immediately see that requiring the family $\cF$ to be $t$-cancellative 
 is equivalent for its representation set of binary vectors to satisfy the following: 
 for every $(t+2)$-tuple $({\mathbf x}^{(1)}, {\mathbf x}^{(2)}, \dots, {\mathbf x}^{(t+2)})$ 
  of distinct vectors in the set (considered in an arbitrary but fixed order) there exist 
  at least $t+1$ different values of $k\in[n]$, 
  such that the corresponding ordered $(t+2)$-tuples $(x_k^{(1)}, x_k^{(2)},\dots, x_k^{(t+2)})$
   are all
 different while for each of them we have the sum $x_k^{(1)} + x_k^{(2)}+\dots +x_k^{(t+2)}= 1$.
In hypergraph language, at least $t+1$ of the sets among the $t+2$ are having degree 
 one vertices.
This problem can be seen in a more general context. 
  We can require that for every ordered $a$-tuple of vectors in the set, 
  there exists at least $b$ different columns, which sum up to 1.
 (obviously $1\leq b\leq a$). 
$$t(a,b):=\limsup_{n\to \infty}\frac{1}{n}\log_2 f(n, \text{locally $(a,b)$-thin}).
 $$
We have $t(a,1)\geq t(a,2)\geq \dots \geq t(a,a)$. 

This problem was investigated by Alon, Fachini, K\"orner, and Monti~\cite{AFK,AKM,FKM},  
  they proved that $t(4,1)<0.4561\dots$ and 
  $t(a,1)< 2/a$ for all even $a$ and 
  $$\Omega(\dfrac{1}{a})\leq t(a,1)\leq O(\dfrac{\log a}{a})< 0.793
  $$ 
for all $a$. 
This is a notoriously hard problem. 
In particular, one does not even know whether $t(3,1) < 1$ 
 (see, Erd\H os and Szemer\'edi~\cite{ESz}). 

Concerning one of the most interesting cases, the case $a=4$, a
  locally $(4,1)$-thin family is also {\em weakly union-free} ($A\cup B=C\cup D$ implies
 $\{A,B \}=\{ C,D\}$). 
The best upper bound 
 $$ \log_2 \, f(n,\,\text{weakly union-free} ) < (0.4998\dots +o(1))\,n$$
 is due to Coppersmith and Shearer~\cite{CS}.  
Nothing nontrivial is known about $t(4,2)$.  
Our Theorem~\ref{th:2c} implies $t(4,3)\leq \log_2 5 -2= 0.3219\dots$
One can find more similar problems in the survey article by K\"orner~\cite{Kor1995}
  and in the more recent paper by K\"orner and Monti~\cite{KM}.

\section{Cancellative and cover-free families}

A family ${\mathcal F}\subseteq 2^{[n]}$ is $g$-{\em cover-free}
  if it is locally $(g+1,g+1)$-thin.
In other words, for arbitrary distinct members $A_0,A_1,\dots , A_g\in {\mathcal F}$
$$A_0\not\subseteq\bigcup_{i=1}^g A_i.$$
Let $C(n,g)$ ($C_r(n,g)$) be the maximum size of a $g$-cover-free
 $n$ vertex code ($r$-uniform hypergraph, resp.).
Obviously, $C(n,g) \leq C(n,g-1)\leq \dots C(n,1)$ 
 and  $C_r(n,g)\leq  C_r(n,g-1)\leq \dots C_r(n,2)$. 

Union free and cover free families were introduced by Kautz and Singleton~\cite{KaSi}. 
They studied binary codes with the property that the disjunctions (bitwise $OR$s) of 
 distinct at most $g$-tuples of codewords are all different. 
In information theory usually these codes are called {\bf superimposed} and they have been
 investigated in several papers on multiple access communication (see,
 e.g., Nguyen Quang A and Zeisel~\cite{AZ}, D'yachkov and Rykov~\cite{DR1},
 Johnson~\cite{J1}).
The same problem has been posed -- in different terms --
 by Erd\H os, Frankl and F\"uredi~\cite{EFF1, EFF2} in
 combinatorics, by S\'os \cite{S} in combinatorial number theory, and
 by Hwang and S\'os~\cite{HS} in group testing. 
For recent generalizations see, e.g.,  Alon and Asodi~\cite{AA1}, and De Bonis
and Vaccaro~\cite{DV}.
D'yachkov and Rykov~\cite{DR1} proved that here are 
 positive constants $\alpha_1$ and $\alpha_2$ such that
\begin{equation}\label{eq:DR}  
    \alpha_1\frac{1}{g^2}< \frac{\log C(n,g)}{n} < \alpha_2\frac{\log g}{g^2}
   \end{equation}  holds for every $g$ and $n> n_0(g)$. 
One can find short proofs of this upper bound in~\cite{F} and in 
 Ruszink\'o~\cite{Ru}.

Using induction on $t$ we extend Theorem~\ref{th:2c} for all $t\geq 2$. 
\begin{thm}\label{th:tcanc} \quad There exists an absolute constant $\alpha >0$
 such that 
$$c(n,t)< \alpha \,  n^{(t-1)/2} \left( \dfrac{t+3}{t+2}\right)^n$$
 holds for all $n,t\geq 2$.  
 \end{thm}
The proof is postponed to Section~\ref{ss:tcanc}. 
This might give a decent upper bound for small $t$ but the true order of 
 magnitude of $c(n,t)$ is much smaller.

\begin{thm}\label{th:t*canc}\quad 
There exist positive constants $\beta_1$ and $\beta_2$ and a bound $n_0(t)$ depending 
 only from $t$ such that the following bounds hold for all $n> n_0(t)$, $t\geq 2$
\begin{equation}\label{eq:t*canc}
   \beta_1\frac{1}{t^2}< \frac{\log c(n,t)}{n} < \beta_2\frac{\log t}{t^2}.
  \end{equation}
 \end{thm}

\noindent
{\em Proof:} \quad
The lower bound is a standard random choice and it follows from (\ref{eq:DR}) since
  $c(n,t)\geq C(n,t+1)$.  
For the upper bound we observe that
\begin{equation}\label{eq:t*}
  c(n,t) \leq  1+ \big\lfloor\frac{t}{2}\big\rfloor + 
      C(n,\big\lfloor\frac{t}{2}\big\rfloor).
  \end{equation}
Indeed, if $\cF\subset 2^{[n]}$ exceeds the right hand side, then one can find
 $h+1$ distinct members $A_0, A_1, \dots, A_h\in \cF$, where $h=\lfloor t/2\rfloor$, 
  such that $A_0\subset A_1\cup \dots \cup A_h$. 
Then, the size of the family $\cF':=\cF\setminus \{ A_0, A_1, \dots, A_h\}$
still exceeds $C(n,h)$ so there is another set of distinct members $B_0, \dots, B_h\in \cF'$
 with $B_0\subset B_1\cup \dots \cup B_h$. 
Taking another set $D\in \cF'$ if $t$ is odd, we have selected $t+2$ distinct members of 
 $\cF$ such that the union of $t$ of them, 
 namely $A_1, \dots, A_h$ and $B_1, \dots, B_h$ and eventually $D$, covers the other two, 
  namely $A_0$ and $B_0$. 
  Hence $\cF$ cannot be $2$-cancellative.
  
Finally, the upper bound (\ref{eq:t*canc}) is implied by (\ref{eq:t*}) and (\ref{eq:DR}). \qed

\section{Three-uniform cancellative families and sparse hypergraphs}

The rest of our results concern about $r$-uniform cancellative families.
We are especially interested in the case when $n$ is large with respect to $r$. 

Frankl and the present author~\cite{FFF} determined asymptotically the maximum 
 size of an $r$-uniform $g$-cover-free family showing that there exists a
 positive constant $\gamma:=\gamma(r,t)$ such that 
\begin{equation}\label{eq:Crg}
   C_r(n,g)= (\gamma + o(1))n^{\lceil r/g\rceil}
   \end{equation}
   as $r$ and $g$ are fixed and $n$ tends to infinity. 
A way to determine $\gamma(r,t)$ was also described. 
This and the $r$-uniform version of (\ref{eq:t*})
\begin{equation*}
  C_r(n,t+1)\leq c_r(n,t) \leq  1+ \big\lfloor\frac{t}{2}\big\rfloor + 
      C_r(n,\big\lfloor\frac{t}{2}\big\rfloor).
  \end{equation*}
imply
\begin{equation}\label{eq:gen_rt}
  (\gamma(r,t+1)-o(1))n ^{\delta_1}\leq c_r(n,t) \leq  
   (\gamma(r,\lfloor t/2\rfloor)+o(1)) n^{\delta_2}.
  \end{equation}
where the exponents are  $\delta_1:=\lceil r/(t+1)\rceil$ and 
$\delta_2:=\lceil r/\lfloor t/2\rfloor\rceil$.
The next theorem shows that to obtain the true asymptotic for $c_r(n,t)$ 
 like the one in (\ref{eq:Crg}) for $C_r(n,g)$ is probably a very difficult problem
 even in the case $r=3$. 

Brown, Erd\H os and  S\'os~\cite{E64, BES1, BES2} introduced the function $f_r(n,v,e)$ to denote
the maximum number of edges in an $r$-uniform hypergraph on $n$ vertices which does not contain
 $e$ edges spanned by $v$ vertices.
Such hypergraphs are called  $\bbG(v,e)$-{\em sparse} (more precisely $\bbG_r(v,e)$-{\em sparse}).
They showed that $f_r(n,e(r-k)+k,e)=\Theta(n^k)$ for every $2 \leq k < r$ and $e \geq 2$. 
The upper bound $(e-1)\binom{n}{k}$ is easy, 
  no $k$-set can be contained in $e$ hyperedges. 
If we forbid $e$ edges spanned by one more vertices then the problem
 becomes much more difficult.
Brown, Erd\H os and  S\'os conjectured that
\begin{equation*} 
   f_r(n,e(r-k)+k+1,e)=o(n^k).
   \end{equation*}
The $(6,3)$-Theorem of Ruzsa and Szemer\'edi~\cite{RSz} deals with the case 
 $(e,k,r)=(3,2,3)$, when no six points contain three triples. 
They showed that there exists an $\alpha> 0$ such that
\begin{equation}\label{eq:63}
     n^2 e^{-\alpha \sqrt{\log n}} =  n^{2-o(1)} <  f_3(n,6,3)=o(n^2).
\end{equation}
Since a $\bbG(6,3)$-sparse system is $\bbG(7,4)$-sparse we have 
\begin{equation}\label{eq:6374}
f_3(n,6,3)\leq f_3(n,7,4)
  \end{equation}
and Erd\H os {\bf conjectures} that $f_3(n,7,4)=o(n^2)$.

\begin{thm}\label{th:2c3}
\begin{equation}\label{eq:10_10}
  f_3(n,7,4) -\frac{2}{5}n\leq c_3(n,2) \leq  \frac{9}{2} f_3(n,7,4) +n. 
  \end{equation}
 \end{thm}
The proof is presented in Section~\ref{ss:r=3}.

The $(6,3)$ thorem was extended by  Erd\H os, Frankl, and R\"odl~\cite{EFR} for arbitrary fixed $r\geq 3$,
\begin{equation}\label{eq:EFR}
 n^{2-o(1)}<f_r(n,3(r-2)+3,3)=o(n^2), 
  \end{equation}
 and then by  Alon and Shapira~\cite{AS}
 $n^{k-o(1)} < f_r(n,3(r-k)+k+1,3) = o(n^k)$.
Even the case $k=2$,  $f_r(n,e(r-2)+3,e)=o(n^2)$,  is still open.
Nearly tight upper bounds were established by  S\'ark\"ozy and
Selkow~\cite{SS2, SSk}:
$
  f_r(n,e(r-k)+k+\lfloor\log_2 e\rfloor,e)=o(n^k)$ for $r>k\geq2\ {\rm and}\ e\geq 3,
  $
and 
$
  f_r(n,4(r-k)+k+1,4)=o(n^k)
  $
for the case $e=4$, $r>k\geq 3$.

\section{An upper bound for uniform families}\label{ss_uniform_upper}

\begin{thm}\label{th:up}
For every $k$ and $n$ we have
\begin{equation*} 
c_{2k}(n,2)\leq \frac{\binom{n}{k}}{\frac{1}{2}\binom{2k}{k}}.
  \end{equation*}
 \end{thm}

To estimate $c_{2k+1}(n,2)$ let us consider a $(2k+1)$-uniform family on $[n]$ 
 and join the element $(n+1)$ to each hyperedge.
If the original family is $t$-cancellative, then so is the extended family. 
We can apply Theorem~\ref{th:up} to get
 \begin{equation}\label{eq:up_odd}
c_{2k+1}(n,2)\leq \frac{\binom{n+1}{k+1}}{\frac{1}{2}\binom{2k+2}{k+1}}.
  \end{equation}

\noindent{\em Proof of Theorem~\ref{th:up}.}\quad 
Suppose that $\cF$ is a $2k$-uniform, 2-cancellative family with the
 underlying set $[n]$.
In case of $|\cF|\leq 3$ the inequality holds, so we may suppose that 
 $|\cF|>3$.
Then $\cF$ is 1-cancellative, too.

Define a graph $G=(V,E)$ with vertex set $V:= \binom{[n]}{k}$, i.e., the 
 family of $k$-subsets of $[n]$.   
A pair $A,B\in V$ forms an edge of $G$ if $A\cup B\in \cF$.
Such a pair necessarily contains disjoint sets.
Since every $F\in \cF$ has $\frac{1}{2}\binom{2k}{k}$ partitions into $k$-sets
 we have
  $$ |\cE(G)|=\frac{1}{2}\binom{2k}{k} |\cF| . $$
We claim that $|\cE(G)|\leq |V|=\binom{n}{k}$.

Consider four adjacent edges in $G$ on five (not necessarily distint) vertices 
 $V_1, \dots, V_5\in V(G)$  (in fact these are $k$-sets of $[n]$)
such that  $\{V_i,V_{i+1}\}\in \cE(G)$ ($1\leq i\leq 4$) and $V_i\neq V_{i+2}$. 
If these four edges determine four distint sets  $V_i\cup V_{i+1}\in \cF$, then 
 the identity
$$
(V_1\cup V_2)\cup (V_4\cup V_5)\cup (V_2\cup V_3)=
  (V_1\cup V_2)\cup (V_4\cup V_5)\cup (V_3\cup V_4)
  $$
 yields a contradiction, since $\cF$ is 2-cancellative. 
By definition we have $(V_1\cup V_2)\neq (V_2\cup V_3)\neq
  (V_3\cup V_4)\neq (V_4\cup V_5)$. 
We also have $V_1\cup V_2\neq V_3\cup V_4$ (and by symmetry $V_2\cup V_3\neq V_4\cup V_5$).
Indeed, $V_3\subset (V_1\cup V_2)$, $V_2\cap V_3=\emptyset$ leads to $V_1=V_3$
 what we have excluded.
The last case to investigate is when $V_1\cup V_2= V_4\cup V_5$ 
  and the four edges determine exactly three sets. 
This leads to the contradiction 
\begin{eqnarray*}
(V_1\cup V_2)\cup (V_2\cup V_3)&=& (V_1\cup V_2) \cup V_3 =(V_4\cup V_5)\cup V_3\\
  &=& (V_4\cup V_5)\cup (V_3\cup V_4)= (V_1\cup V_2)\cup (V_3\cup V_4). 
  \end{eqnarray*}
We conclude that $G$ does not have such a sequence of four edges.
Therefore $G$ contains no cycles, neither paths of length 4, implying
 $|\cE(G)| < |V|$. \qed

\section{The nonuniform case, the proof of  Theorem~\ref{th:2c}}\label{ss:2c}

Suppose that $\cF$ is a 2-cancellative family of maximal size on the
 underlying set $[n]$.
Split $\cF$ according the sizes of its edges
$\cF_r:= \{ F\in \cF: |F|=r\}$. 

The sequence $\sqrt{2k-1}\binom{2k}{k}4^{-k}$ is monotone increasing for $k=1,2,3\dots$
 so we obtain that 
 $\binom{2k}{k}^{-1}\leq 2\sqrt{2k-1}\times 4^{-k}$ for all $k\geq 1$. 
Using this upper bound in  Theorem~\ref{th:up} 
  we obtain
$$
  c_{2k}(n,2)\leq \frac{\binom{n}{k}}{\frac{1}{2}\binom{2k}{k}}
\leq \binom{n}{k} 4^{-k} 4 \sqrt{2k-1}<  4\sqrt{n}\times  
  \binom{n}{k} 4^{-k}. 
  $$
The same inequalty and (\ref{eq:up_odd}) give
$$
  c_{2k+1}(n,2)\leq \frac{\binom{n+1}{k+1}}{\frac{1}{2}\binom{2k+2}{k+1}}
\leq \binom{n+1}{k+1} 4^{-k-1} 4 \sqrt{2k+1}\leq  4\sqrt{n}\times  
  \binom{n+1}{k+1} 4^{-k-1}. 
  $$
Finally,
\begin{eqnarray*}
c(n,2)=|\cF|&=&\sum_r |\cF_r|\leq \sum_r c_r(n,2)= \left( \sum_{k\geq 0} c_{2k}(n,2)\right)
   +  \left( \sum_{k\geq 0} c_{2k+1}(n,2)\right)\\
&<& \left( \sum_{k\geq 0} 4\sqrt{n}\times  
  \binom{n}{k} 4^{-k}\right) + \left( \sum_{k\geq 0} 4\sqrt{n}\times  
  \binom{n+1}{k+1} 4^{-k-1}\right) \\
&=& 4{\sqrt{n}}\left( \left(1+\frac{1}{4}\right)^n+\left( 1+\frac{1}{4} \right)^{n+1} \right) = 9\sqrt{n} \left(\frac{5}{4}\right)^n.
 \end{eqnarray*}

\section{The case of $t$-cancellative codes, the proof of 
  Theorem~\ref{th:tcanc}}\label{ss:tcanc}

We will define a monotone sequence $0< \alpha_2\leq \alpha_3\leq \dots \alpha_t\leq \dots $
 which is bounded above (by $\alpha$) such that 
  \begin{equation}\label{eq:tc}
   c_t(n,2)<   \alpha_t \, n^{(t-1)/2} \left( \frac{t+3}{t+2}\right)^n
  \end{equation}
holds for every $n,t\geq 2$. 
By Theorem~\ref{th:2c} this holds for $t=2$ with $\alpha_2:=9$. 
Suppose that $t\geq 3$ and (\ref{eq:tc}) holds for $t-1$. 
We use the following upper bound
\begin{equation}\label{eq:rough}
  c_r(n,t)\leq c (n-r, t-1).
  \end{equation}
Indeed, if $\cF$ is an $r$-uniform $t$-cancellative family on $[n]$,
 then for any $F_0\in \cF$ the family 
  $\{ F\setminus F_0: F\in \cF, F\neq F_0\}$
 is $(t-1)$-cancellative. 
Use the inequality
$$
  \binom{r(t+2)}{r} > \frac{1}{3\sqrt{r}} \left(\frac{(t+2)^{t+2}}{(t+1)^{t+1}} \right)^{r}
  $$
which holds for every integers $r\geq 1$, $t\geq 0$ and substitute 
  $n=r(t+2)$ to (\ref{eq:rough}). 
We obtain
$$
  \frac{c_r(r(t+2),t)}{\binom{r(t+2)}{r}}< \alpha_{t-1}\, n^{(t-2)/2} \left( 
    \frac{t+2}{t+1}\right)^{r(t+1)}\times 3\sqrt{r} \frac{(t+1)^{r(t+1)}}{(t+2)^{r(t+2)}}
    = \frac{3\alpha_{t-1}}{\sqrt{t+2}}  n^{(t-1)/2} (t+2)^{-r}.
  $$
For $n\geq m$ we have $c_r(n,t)\binom{n}{r}^{-1}\leq c_r(m,t)\binom{m}{r}^{-1}$.
We obtain that the right hand side is an upper bound for 
$c_r(n)/\binom{n}{r}$ for every $n\geq r(t+2)$.
 For any given $n$ and $t$ this gives  
$$
  \sum_{r\leq n/(t+2)} c_r(n,t)\leq \sum \binom{n}{r} \frac{3\alpha_{t-1}}{\sqrt{t+2}}  n^{(t-1)/2} (t+2)^{-r} \leq  \frac{3\alpha_{t-1}}{\sqrt{t+2}}  n^{(t-1)/2} \left(1+\frac{1}{t+2}\right)^{n}. 
  $$
We estimate the case $n< r(t+2)$ using (\ref{eq:rough}) again
\begin{eqnarray*}
   \sum_{n/(t+2)< r\leq n} c_r(n,t)&\leq& \sum_{r> n/(t+2)} c (n-r,t-1)\enskip
     < \sum_{n-r< (t+1)n/(t+2)}  \alpha_{t-1} \, n^{(t-2)/2} \left( \frac{t+2}{t+1}\right)^{n-r}
 \\ &<&  \alpha_{t-1} \, n^{(t-2)/2} (t+2) \left( \frac{t+2}{t+1}\right)^{n(t+1)/(t+2)}.
  \end{eqnarray*}
Here $\left( \frac{t+2}{t+1}\right)^{(t+1)/(t+2)}< (t+3)/(t+2)$ so the sum of the above two 
 displayed formula gives
$$
 c(n,t)\leq \alpha_{t-1}\left(   \frac{3}{\sqrt{t+2}}  +  \frac{t+2}{\sqrt{n}} \right) 
   n^{(t-1)/2} \left( \frac{t+3}{t+2}\right)^{n}.
  $$
The rest is a little calculation (e.g., we may suppose that $n> 2(t+2)^2$, otherwise
 our upper bound (\ref{eq:tc}) for $c(n,t)$ exceeds the same upper bound for $c(n,t-1)$.)
\qed

\section{Three-partite hypergraphs}\label{ss:r=3}

In this section we prove Theorem~\ref{th:2c3} on 3-uniform 2-cancellative families.

\begin{lemma}[Erd\H os and Kleitman~\cite{EK}]\label{le:EK}\quad 
Let $\cF$ be an $r$-uniform hypergraph.
Then there exists an $r$-partite $\cF^*\subset \cF$ with 
$ |\cF^*|\geq \dfrac{r!}{r^r}|\cF|$.  \qed
  \end{lemma}

\noindent
{\em Proof of Theorem \ref{th:2c3}:} \quad
Suppose that $\cH$ is a 3-uniform $\bbG(7,4)$-sparse family with vertex set $[n]$.  
We claim that there exists a subfamily $\cH''\subset \cH$ 
 such that 
$$
   |\cH''|\geq |\cH|-\frac{2}{5}n \quad \text{and}\quad \cH''\text{ is linear}. 
  $$
First, note that if the hyperedge $F\in \cH$ has 
 two other edges $F_1, F_2$ with $|F\cap F_i|=2$, then these three edges
 form a separate connected component of $\cH$ on $5$ vertices.
Let $\cH'$ be the hypergraph obtained from $\cH$ after deleting 
 two edges from each such 5-vertex component.
If $F_1\in \cH'$ and there exists an $F_2\in \cH'$ with $|F_1\cap F_2|=2$
 then this $F_2$ is unique. 
Moreover, if $|F_1\cap F_2|=2$ and $|F_3\cap F_4|=2$ hold for four distint 
 sets, then $F_1\cup F_2$ is disjoint to $F_3\cup F_4$.
Leave out an edge from each such a pair to obtain $\cH''$, which is 
 a linear hypergraph by definition, and we left out at most $(2/5)n$ edges of $\cH$. 

Second, observe that a linear, $3$-uniform, $\bbG(7,4)$-sparse family $\cH''$ is 
 $2$-cancellative.
Indeed, if $X:=C\setminus (A\cup B)= D\setminus (A\cup B)$ for 
 four distinct members $\{A,B,C,D\} \subset\cH''$, then $X\subset C\cap D$ so $|X|\leq 1$,
 and $|A\cup B|\leq 6$, so they form a $\bbG(7,4)$ family, a contradiction.
If we take $|\cH|$ as large as possible, then we complete the proof of 
 the first inequality of (\ref{eq:10_10}) as follows
$$
  c_3(n,2)\geq |\cH''|\geq |\cH|-\frac{2}{5}n = f_3(n,7,4)-\frac{2}{5}n.
  $$

Next, let $\cF$ be a 3-uniform, 2-cancellative family on $n$ vertices.
We claim that there exists a subfamily $\cF'\subset \cF$ 
 such that 
$$
   |\cF'|\geq |\cF|-n \quad \text{and}\quad \cF'\text{ is linear}. 
  $$
Indeed, leave out a hyperedge from $\cF$ if it has a vertex of degree one.
Repeat this until we get $\cF'\subset \cF$ with every degree is either 0 or at least 2.
We claim that $\cF'$ is linear (in case of $|\cF'|\geq 4$). 
Suppose not, $|F_1\cap F_2|=2$, $F_1, F_2\in \cF'$, 
 $x_i$ is the unique element of $F_i\setminus(F_1\cap F_2)$. 
By our degree condition there exist $A_i\in \cF'$, $x_i\in A_i$, $A_i\neq F_1$ and
$A_i\neq F_2$.
This leads to the contradiction $A_1\cup A_2\cup F_1=A_1\cup A_2\cup F_2$.
(The case $|\cF'|=3$ is left to the reader.) 

Apply Lemma~\ref{le:EK} to $\cF'$ to obtain a $3$-partite $\cF^*$ of size 
 $|\cF^*|\geq (2/9)|\cF'|$. 
We claim it is $\bbG(7,4)$-sparse, because
 every $3$-partite, $2$-cancellative, linear family $\cF^*$ is
 $\bbG(7,4)$-sparse. 
Indeed, take any four distict members $\{A,B,C,D\} \subset\cF^*$.
If $A\cap B=\emptyset$ then $C\setminus (A\cup B)\neq \emptyset$ (by linearity)
and it is not equal to $D\setminus (A\cup B)$ (by 2-cancellativeness)
 so the union of the four of them have at least 8 vertices.
Otherwise they pairwise have a one-element intersection.
If there is a degree three, say $A\cap B\cap C=\{ x\}$, then linearity and 
 3-partiteness imply that $D$ is not covered by $A\cup B\cup C$, so
 again their union has eight vertices (at least).  
If these four triples pairwise meet but their maximum degree
 is two, then 
$C$ and $D$ covers $(A\cup B)\setminus (A\cap B)$, and they have a 
 single common vertex outside $A\cup B$ yielding the contradiction
 $A\cup B\cup C=A\cup B\cup D$.

Finally, if we take $|\cF|$ as large as possible, then we complete the proof of 
 the second inequality of (\ref{eq:10_10}) as follows
$$
  f_3(n,7,4)\geq |\cF^*|\geq \frac{2}{9}|\cF'|\geq \frac{2}{9}(|\cF|-n)
    = \frac{2}{9}(c_3(n,2)-n).  \qed
  $$

Define the hypergraphs $\bbG_6$  and $\bbG_7$
 as follows on $6$ and $7$ vertices.\\ \indent\indent
$\cE(\bbG_6):=\{ 123, 156, 426, 453\}$, \\ \indent\indent
$\cE(\bbG_7):=\{ 123, 456, 726, 753\}$.\\
Note that both are three-partite and the 3-partition of their vertices is unique.
We have proved the following fact.

\begin{proposition}\label{le:3-partite_2canc} \quad
Suppose that $\cF$ is a three-partite, linear hypergraph. 
It is $2$-cancellative if and only if it avoids $\bbG_6$, and $\bbG_7$.
It is $\bbG(7,4)$-sparse if and only if it avoids $\bbG_6$, and $\bbG_7$. \qed
  \end{proposition}

\section{A construction by induced packings}

According to the upper bounds in Theorem~\ref{th:up} 
we have
$$
 c_2(n,2)\leq n,\enskip\quad c_3(n,2)\leq \frac{1}{6}n(n+1),\enskip\quad c_4(n,2)\leq \frac{1}{6}n(n-1).
  $$
Obviously, $c_2(n,2)=n-1$ for $n>3$. 
The second inequality, although it is close to the true order of magnitude, is not sharp
 if Erd\H os conjecture is true, see (\ref{eq:63}), (\ref{eq:6374}), and (\ref{eq:10_10}).
Any 4-uniform Steiner system $S(n,4,2)$ is 2-cancellative yielding the lower bound
 $c_4(n,2)\geq \frac{1}{12}n(n-1)-O(n)$ for all $n$.

\begin{thm}\label{th:2c4}\quad
$c_4(n,2)=\dfrac{1}{6}n^2-o(n^2)$.
 \end{thm}

\begin{thm}\label{th:2ck}\quad
$c_{2k}(n,k)\geq \dfrac{n^k}{k^k}-o(n^k). $
 \end{thm}

The proof of Theorem~\ref{th:2ck} is postponed to Section~\ref{ss:proof_of_2k}, 
 for the construction giving Theorem~\ref{th:2c4} we use induced packings of graphs.

A set of graphs $\cP:=\{G_1=(V_1, {\mathcal E}_1)$, $G_2=(V_2, {\mathcal E}_2)\dots\}$
 is called a {\em packing} if they are edge-disjoint subgraphs of $G=(V, {\mathcal E})$
  (by definition $V_i\subset V$ for each $i$).
The packing $\cP$ is called an {\em induced packing} if $G$ restricted to $V_i$ is 
 exactly $G_i$ (for all $i$).
The induced packing $\cP$ is called an {\bf almost disjoint induced packing} into 
 the graph  $G$ if $|V_i\cap V_j|\leq 2$ (for all $i\neq j$).
It follows that whenever $F=V_i\cap V_j$, $|F|=2$, then $F$ is not an edge of $G$. 
Shortly, any two induced graphs $G[V_i]$ and $G[V_j]$ are either vertex disjoint, or 
 share one vertex, or meet in a non-edge. 
For example, if $V$ is an $n$-element set, $n$ is even, $V=A_1\cup A_2\cup \dots \cup A_{n/2}$ 
 where each $|A_i|=2$ and $G$ is the complete graph on $V$ minus the $n/2$
 edges of the perfect matching $\{ A_1, A_2, \dots\}$, then ${\mathcal E}(G)$ 
 can be decomposed into $n(n-2)/8$ almost disjoint, induced four-cycles,
 namely those induced by $A_i\cup A_j$. 
 
Let $H$ be a graph of $e$ edges and let $i(n,H)$ denote the maximum number of almost 
 disjoint induced copies of $H$ can be packed into any $n$-vertex graph.
It was proved by Frankl and the present author that
\begin{equation*}  
  i(n,H)= \frac{1}{e(H)} \binom{n}{2} -o(n^2).
  \end{equation*}
In other words
\begin{lemma} \label{le:FF} {\rm \cite{FFF}}\quad
For any fixed graph $H$ with $e$ edges one can delete
 $o(n^2)$ edges of the graph $K_n$ such that the rest of the edges, the graph
 $L_n=L_n(H)$,  can be decomposed into $(1-o(1)){n \choose 2}/e$ 
 almost disjoint induced copies of $H$. \qed
\end{lemma}

\noindent{\em Proof of Theorem~\ref{th:2c4}.}\quad
The graph $H_k$, for $k\geq 3$, is defined as a complete graph $K_k$ and  
 $\binom{k}{3}$ vertices of degree three, each of those connected to a different 
 triple of $V(K_k)$.
We have that $H_k$ has $k+\binom{k}{3}$ vertices, $\binom{k}{2}+ 3\binom{k}{3}$ edges, 
 and it contains $\binom{k}{3}$ {\em special} $K_4$'s, those having a vertex of degree $3$.
Take any almost disjoint packing of copies of $H_k$, $\cP:=\{ H_k^1,  H_k^2, \dots \}$ and define
 a 4-uniform family $\cF(\cP)$ as the vertex sets of the special $K_4$'s.
Obviously $|\cF(\cP)|= \binom{k}{3}|\cP|$.

We claim that $\cF(\cP)$ is a 2-cancellative family. 

Suppose, on the contrary, that there are four distinct members $A,B,C,D\in \cF$
 with $A\cup B\cup C=A\cup B\cup D$. 
Note that $|F\cap F'|\leq 2$ for $F,F'\in \cF$. Furthermore, in the case of equality
 $F$ and $F'$ are generated by the same $H_k^i$. 
 
Consider first the case when $C$ and $D$ are generated by the same copy of $H_k^i$.
$C=\{ c,x_1,x_2,x_3\}$ $D=\{ d,y_1,y_2,y_3\}$ where $c$ and $d$ are distinct 
  $3$-degree vertices of $H_k^i$. 
The element $c$ is covered by $A\cup B$, say $c\in A$. 
By definition, the pairs $cx_1$, $cx_2$, $cx_3$ are only covered by 
  $C$ among the members of $\cF$ and since 
  there is no edge (of any $H_k^j$) from $c$ to $D\setminus C$ 
  those pairs are not covered by any member of $F\in \cF$.
Hence $A\cap(C\cup D)=\{ c\}$.
Similarly, $B\cap (C\cup D)=\{ d\}$.
Hence $(A\cup B)\cap (C\cup D)=\{ c,d\}$. 
Since the symmetric difference $C\Delta D$ is contained in $A\cup B$, we obtain that
 it is $\{ c,d\}$. 
This leads to the contradiction $|C\cap D|=3$. 

Consider now the other case, that $i\neq j$ and $C=\{ c,x_1,x_2,x_3\}$ is generated 
 by $H_k^i$ where $c$ is a $3$-degree vertex of $H_k^i$, and 
 $D=\{ d,y_1,y_2,y_3\}$ is generated by $H_k^j$ where $d$ is a
  $3$-degree vertex of $H_k^j$. 
Since they come from different copies of $H_k$ we have $|C\cap D|\leq 1$.
This implies that either $A$ or $B$ meets $C$ in two vertices, say $|A\cap C|=2$.
Then $A$ is generated by $H_k^i$ as $C$ is, say $A=\{ a,x_2,x_3,x_4\}$.  
It follows that $|A\cap D|\leq 1$ so $|B\cap D|=2$ and $D$ is generated by $H_k^j$ as well. 
The pair $\{ x_1, c\} =C\setminus A$ is covered by $B\cup D$, thus it is
 covered by $V(H_k^j)$. 
This leads to the contradiction that $V(H_k^i)\cap V(H_k^j)$ contains an edge
 of $H_k^i$. 
This completes the proof that $\cF$ is 2-cancellative.

For given $n$ taking a large induced packing of $H_k$'s  Lemma~\ref{le:FF} implies
 that 
\begin{equation}\label{eq:c4_Hk} 
   c_4(n,2)\geq i(n,H_k)\binom{k}{3}\geq (1-o(1))
    \frac{\binom{k}{3}}{\binom{k}{2}+ 3\binom{k}{3}}\binom{n}{2}
    \end{equation}
when $k$ is fixed and $n\to \infty$.
Let $\pi_4:=\liminf_{n\to \infty}\{ c_4(n,2)/\binom{n}{2}\}$. 
The lower bound (\ref{eq:c4_Hk}) gives that 
 $$  \pi_4\geq  \frac{1}{\binom{k}{2}/\binom{k}{3} + 3}.$$
Since this holds for each $k$ we obtain $\pi_4\geq 1/3$.
Finally, $\pi_4\leq 1/3$ was proved in Theorem~\ref{th:up}, completing 
 the proof of $\pi_4=1/3$. \qed

\section{The lower bound for the $2k$-uniform case}\label{ss:poly}

\subsection{Not vanishing polynomials}

For a set of variables $X=\{ x_1, \dots, x_s\}$ and $0\leq i\leq s$
 the symmetric polynomial $\sigma_i(X)$ is defined as
  $\sum_{I\subset X, |I|=i}\prod _{i\in I}x_i$, 
 $\sigma_0(X)=1$. 
For convenience, frequently $\sigma_i(X)$ is defined to be $0$ for
  $|X|< i$ (and for $i<0$).
Suppose that $X_1$, $X_2, \dots, X_\ell$ are disjoint sets of variables,
 $|X_j|=t_j$, $0<t_j< k$, $\sum_j (k-t_j)=k$.
The entries of a row of the $k\times k$ matrix $M(X_1, \dots, X_\ell)$
 consists of a block with the symmetric polynomials  
  $\{\sigma_0(X_j), \sigma_1(X_j),  \dots, \sigma_{t_j}(X_j)\}$ and 
  $0$'s otherwise.
The rows are distinct, so these blocks are shifted in all possible 
 $k-t_j$ ways.
\begin{equation*}
M(X_1, \dots, X_\ell):= \left(
\begin{array}{ccccccccc} 
1& \sigma_1(X_1)&\sigma_2(X_1) &\dots &\dots  &\sigma_{t_1}(X_1) & 0& \dots & 0\\
0 &1& \sigma_1(X_1)&\sigma_2(X_1) &\dots &\dots  &\sigma_{t_1}(X_1) & 0& \vdots \\ 
{}\\
\vdots& \vdots &\ddots & \ddots& & & &\ddots &0 \\
0 &0 &0 &1& \sigma_1(X_1)&\sigma_2(X_1) &\dots &\dots  &\sigma_{t_1}(X_1)\\
1& \sigma_1(X_2)&\sigma_2(X_2) &\dots &\dots  &\sigma_{t_2}(X_2) & 0& \dots & 0\\
0& \ddots &\ddots & \ddots& & & &\ddots &0 \\
{}\\& & & & & & & & 0\\
0 &0 &0 &1& \sigma_1(X_2)&\sigma_2(X_2) &\dots &\dots  &\sigma_{t_2}(X_2)\\
\vdots& \vdots & & & & & & &\vdots  \\
{}\\
1& \sigma_1(X_\ell)&\sigma_2(X_\ell) &\cdots & & \cdots &\sigma_{t_\ell}(X_\ell) & \cdots & 0\\ 
\vdots& \ddots &\ddots & \ddots& & & &\ddots &0 \\
& & & & & & & & \\
0 &\cdots &1 & \sigma_1(X_\ell)&\sigma_2(X_\ell)& & &\cdots   &\sigma_{t_\ell}(X_\ell)\\
\end{array}
\right)
\end{equation*}

\begin{fact}\label{fa:1}\quad
  The polynomial $\det M(X_1, \dots, X_\ell)$ of $\sum |X_j|$ variables is not vanishing. 
  \end{fact}

\noindent{\em Proof.}\quad
Over any field we can substitute only $1$'s and $0$'s such that the matrix $M$
 becomes a lower triangular matrix, having only 1's in the main diagonal and 0's above. 
Namely, let $x=0$ for each $x\in X_1$.
In the second block of $M$, in rows $(k-t_1)+1$ to $(k-t_1)+(k-t_2)$
 only  $\sigma_i(X_2)$ stands in the main diagonal, where $i=k-t_1$.
Define $k-t_1$ variables of $X_2$ to be $1$, the rest $0$.

In general, in the $j$'th block, in rows $(k-t_1)+ \dots +(k-t_{j-1})+1$
 to $\sum_{1\leq s\leq j} (k-t_s)$
 we define the variables of $X_j$ such a way that 
 $(k-t_1)+ \dots +(k-t_{j-1})$ of them are 1's and the rest are 0's.
This can be done, since $(k-t_1)+ \dots +(k-t_{j-1}) \leq t_j$. 
\qed

One can define the matrix in a more general setting
when the blocks consist of rows of the form $(\sigma_{m+1}(X_j), \sigma_{m+2}(X_j),
  \dots, \sigma_{m+k}(X_j))$.
We can get, e.g., that the determinant of the $k\times k$ matrix
 $M$ with $M_{i,j}:= \sigma_{m+i+j-2}(X)$ is not vanishing 
if $m\geq 0$, $|X|\geq m+k-1$. 

Let $q>1$ be a power of a prime, $\bF:=\bF_q$ the finite field of size $q$.
For any polynomial $p(x_1, \dots, x_s)$ over this field the zero set $Z(p)$
 is defined as $Z:= \{(x_1, \dots, x_s)\in {\mathbf F}_q^s: p(x_1, \dots, x_s)=0\}$. 
It is well-known that

\begin{fact}\label{fa:2}\quad
There is a $\eta:=\eta(h,s)$ such that there is an upper bound for the size
 of the zero-set in ${\mathbf F}_q^s$ of any $s$ variables, degree $h$, 
 not identically $0$ polynomial $p(x_1, \dots, x_s)$ as
 $$
  | Z(p(x_1, \dots, x_s))|\leq \eta(h,s) q^{s-1}.   \qed
    $$
  \end{fact}

\subsection{A lemma on independent polynomials}\label{ss:indep}

Let $k$ be a positive integer, let 
$\cP:=\cP_{<k}[{\bf F}, x]$ be the ring of polynomials of degree at most $k-1$
$$
  \cP_{<k}:= \{ a_0+ a_1x+ \dots+ a_{k-1}x^{k-1}: a_i\in \bF\}.
 $$
The number of such polynomials is $q^{k}$ and they form a linear space of 
 dimension $k$ over $\bF$.
A set of polynomials $p_1(x), \dots, p_\ell(x)\in \cP$
 is called $(k_1, \dots, k_\ell)$-{\em independent}, where $k_1, \dots, k_\ell$
 are positive integers if 
 $$
   f_1(x)p_1(x)+ \dots + f_\ell(x)p_\ell(x)\equiv 0
    $$  
    and $\deg (f_i)< k_i$ for all $i$ imply that each $f_i(x)$ is the $0$ polynomial.
Equivalently, all the $q^{\sum k_i}$ polynomials of the form 
  $\sum f_ip_i$ (with $\deg(f_i)< k_i$) are distinct. 
The case all $k_i=1$ corresponds to the usual linear independentness.
To stay in the space $\cP_{<k}$ we also suppose that $\deg(p_i)+k_i< k$.    
Then necessarily $\sum_i {k_i}\leq k$. 

For $Z\subset \bF$ there is a unique polynomial 
 with leading coefficient 1 and roots $Z$, namely
 $$
    p_Z(x):=
                \prod_{z\in Z} (x-z). 
    $$
Suppose that $\ell\geq 2$, $k_1, \dots , k_\ell$ are positive integers with
 $k_1+\dots +k_\ell=k$
 and let $x_1, \dots, x_{(\ell-1)k}$ a sequence of elements of $\bF_q$.
Define the (multi)sets $X_i$ of size $k-k_i$ as intervals of this sequence,
$X_1:=\{ x_s: 1\leq s\leq k-k_1 \}$, in general 
$X_j:=\{ x_s:   \sum_{i< j} (k-k_i)    < s\leq \sum_{i\leq j} (k-k_i)\}$. 
 
\begin{lemma}\label{le:indep}\quad 
The polynomials  $p_{X_1}(x)$, \dots, $p_{X_\ell}(x)$ are $(k_1, \dots, k_\ell)$ 
 independent for all but at most 
$$\left(\eta\left(\binom{k}{2},(\ell-1)k\right)
   +\binom{(\ell-1)k}{2}\right) \, q^{(\ell-1)k-1}
 $$ 
sequences. 
   \end{lemma}
\noindent{\em Proof.}\quad 
We only consider the sequences with distinct elements.
There are at most $\binom{(\ell-1)k}{2} q^{(\ell-1)k-1}$
 sequences with repeated entries.
 
The polynomials $p_1, \dots, p_\ell$ are $(k_1, \dots, k_\ell)$ 
 independent if and only if the set of polynomials 
 $\{x^ip_j(x): 0\leq j< k_j \}$ are linearly independent.
These $k$ polynomials are linearly independent if and only if 
 their coefficient matrix is nonsingular.
The coefficients of $p_{X_j}$ are the values of the symmetric polynomials 
 $\sigma_s(-X_j)$, so the coefficient matrix is exactly
 $M(-X_1, \dots, -X_\ell)$ defined in the previous subsection.
According to Fact~\ref{fa:1} its determinant is a nonzero polynomial
 of degree at most $0+1+2+\dots +(k-1)$ and with  $(\ell-1)k$ variables.
Then Fact~\ref{fa:2} gives an upper bound for the number of 
 sequences whose determinants vanish.  \qed
 
\begin{corollary}\label{co:3indep}\quad
For every $k$ there exists a $q_0(k)$ such that if $q> q_0(k)$ then 
 there exists a $2k$-element set $S\subset \bF_q$ such that 
 the polynomials  
$$p_X(x) , \, p_Y(x),\,   p_W(x) \text{ are }(k-|X|, k-|Y|, k-|W|)
 \text{  independent} $$
 for every partition of $S=X\cup Y\cup W$, $|X|+|Y|+|W|=2k$, $1\leq |X|, |Y|, |W|, < k$. 
    \end{corollary}
    
In fact, applying the previous Lemma with $\ell=3$ we can see that 
 almost all $2k$-sets, all but at most $O(q^{2k-1})$ of them, have this 
 total independentness property. \qed

\subsection{The algebraic construction yielding Theorem~\ref{th:2ck}}\label{ss:proof_of_2k}

Let $q$ be the largest prime power  not exceeding  $n/(2k)$.
Since there are no too large gaps among primes we have $q> n/(2k)- O(n^{5/8})$.
We also suppose that $q> q_0(k)$ used in Corollary~\ref{co:3indep}. 
We are going to define a $2$-cancellative, $2k$-uniform family $\cF$ of size $q^k$.

Take a set $S\subset \bF_q$ of size $2k$ satisfying  Corollary~\ref{co:3indep}.
Our hypergraph $\cF:= \cF(q,S)$ consists of the graphs of the polynomials $\cP_{<k}$
 restricted to $S$.
$V(\cF):=S\times \bF = \{ (s,y): s\in S, y\in \bF \}$, 
every $p\in \cP$ defines a set $F(p):= \{ (s,p(s)): s\in S\}$
and let  $\cF:= \{ F(p): p\in \cP\}$.

To show that $\cF$ is 2-cancelative suppose, on the contrary, that
 $A,B,C$, and $D$ are four distinct members of $\cF$ with 
 $A\cup B\cup C=A\cup B\cup D$.
There are four distinct polynomials $a(x), b(x), c(x)$, and $d(x)\in \cP$
 generating these sets, $A=F(a)$, $B=F(b)$, etc.

Let $W\subset S$ be the set of coordinates where $C$ and $D$ meet,
$W:=\{ s\in S: c(s)=d(s)\}$.
Let $X:=\{ x\in S\setminus W: c(x)=a(x)\}$, and let $Y:=S\setminus (X\cup W)$.
For $x\in X$  $(x,d(x))$ is not covered by $C$ neither $A$, so it must 
 belong to $B$, $b(x)=d(x)$.
For $y\in Y$ we have $c(y)\neq d(y)$, $c(y)\neq a(y)$ so  $(y,c(y))$ must be in $B$,
   $c(y)=b(y)$.
Considering the same $y\in Y$  the element $(y,d(y))$ is not covered by $C$
 neither $B$ so it must belong to $A$, $a(y)=d(y)$.
Let us summarize: there exists a partition of $S=W\cup X\cup Y$ such that
\begin{eqnarray}
  c(w)&=&d(w)\,\text{ for }\, w\in W, \label{eq:111}\\
   c(x)&\neq &d(x)\,\text{ for }\, x\in X,\,\text{ but }\, c(x)=a(x)\,\text{ and }\,d(x)=b(x), \label{eq:112}\\
c(y)& \neq &d(y)\,\text{ for }\, y\in Y, \,\,\text{ but }\,d(y)=a(y) \,\,\text{ and }\,\, c(y)=b(y). 
\label{eq:113}
  \end{eqnarray} 
Since $c$ and $d$ are distinct polynomials of degree at most $k-1$ we have 
 $|W|<k$. 
Similarly $|X|, |Y|< k$.
These also imply that $|X|, |Y|, |W|\geq 2$ (and thus $k\geq 3$).

By (\ref{eq:111}) $c-d$ is divisible by $p_W$, there exists a
 polynomial $c_1(x)\in \cP$ such that
\begin{equation*} 
    c\, =\, d+ c_1p_W\,\text{ where }\, c_1\in \cP,\,\text{ and }\, \deg(c_1)< k-|W|. 
    \end{equation*}
The first halves of (\ref{eq:112}) and (\ref{eq:113})
 similarly imply that 
\begin{eqnarray*}
   a&=&c+ a_1p_X\,\text{ where }\, a_1\in \cP,\,\text{ and }\, \deg(a_1)< k-|X|, 
   \\
   d&=&a+ a_2p_Y\,\text{ where }\, a_2\in \cP,\,\text{ and }\, \deg(a_2)< k-|Y|. 
  \end{eqnarray*} 
Adding these three equations we obtain
\begin{equation*} 
   0=c_1p_W + a_1 p_X + a_2 p_Y. 
    \end{equation*}
Then the independentness of $p_X$, $p_Y$ and $p_W$ imply $c_1=a_1=a_2=0$, a contradiction. \qed

One might think that if we use the second halves of (\ref{eq:112}) and (\ref{eq:113}) 
 then we have more constraints and maybe we do not really need independentness 
 and Corollary~\ref{co:3indep}.
In fact, independentness is essential. 
The second halves only imply that $b= d-a_1p_X = c+a_2p_Y$, so
$\cF$ can have many non-2-cancellative fourtuples if $S$ is not chosen properly.

\section{A remark on $1$-cancellative uniform families}

An $r$-partite hypergraph is cancellative, it contains no three distinct edges with 
$A\cup B=A\cup C$. 
Considering the complete $r$-partite hypergraph on $n$ vertices with almost 
 equal parts one gets
\begin{equation} \label{eq:Boll_conj}
  c_r(n)\geq \lfloor\frac{n}{r}\rfloor\times \lfloor\frac{n+1}{r}\rfloor
    \times \dots \times\lfloor\frac{n+r-1}{r}\rfloor=:p_r(n).
  \end{equation}
The right hand side is exactly $n^r/r^r$ when $r$ divides $n$.
An old result of Mantel, about the maximum size of triangle-free graphs,
 gives $c_2(n)=p_2(n)=\lfloor n^2/4\rfloor$.
Katona~\cite{Kat} conjectured and Bollob\'as~\cite{Boll} proved 
 that $c_3(n)=p_3(n)$.
Bollob\'as also conjectured that 
 equality holds in (\ref{eq:Boll_conj}) for all $n\geq r\geq 4$ as well.
This was established for $2r\geq n\geq r$ in~\cite{Fur27}.
Sidorenko~\cite{Sid4} proved Bollob\'as conjecture for $r=4$. 
(There is a recent refinement of this by Pikhurko~\cite{Pik}).
However, Shearer~\cite{Sh} gave a counterexample. 
His result implies that there exist an $\varepsilon >0$ and $n_0(r)$ 
 such that $c_r(n)> (1+\varepsilon)^r(n/r)^r$ for $n> n_0(r)$, $r\geq 11$.
The cases $5\leq r\leq 10$ are still undecided.

It was observed in~\cite{Fur27}
 that $c_r(n)=2^{n-r}$ for $2r\geq n\geq r$ and
 if $\mathcal{F}$ is a cancellative family of $r$-sets from an $n$-set and $n\geq 2r$, then 
$$ |\mathcal{F}|\leq \frac{2^r}{\binom{2r}{r}}\binom{n}{r}. $$
Here we show an almost matching lower bound

\begin{thm}\label{th:cr_unif}
For every $n\geq r\geq 2$
$$ c_r(n) > \frac{c_0}{2^r}{n \choose r}
  $$
where $c_0:= \prod_{k\geq 1} \frac{2^k-1}{2^k}={.2887\dots}$
 \end{thm}
This result is basically due to Tolhuizen~\cite{T}, although he was not interested
 in $r$-uniform hypergraphs and wrote it like
``the rate of a cancellative code is $\frac{\log 3}{\log 2}-1=.5849\dots$''.
His publication is not even reviewed in MathScinet so we briefly 
 describe his construction.

\noindent{\em Proof.}\quad 
If  $M$ is a random $m\times m$ matrix with entries from the two-element field 
 ${\mathbf F}_2=\{0,1 \}$, then 
$$
  {\text{Prob}(M\text{ is nonsingular})=}\frac{2^m-1}{2^m}\times\frac{2^m-2}{2^m}
   \times \dots \frac{2^m-2^{m-2}}{2^m}\times\frac{2^m-2^{m-1}}{2^m}>c_0.
  $$ 
Considering $(n-r)\times n$ random matrices we obtain 
 an  $(n-r)\times n$ matrix $A$ (over ${\mathbf F}_2$) containing at least 
 $ c_0{n\choose n-r}$ nonsingular $(n-r)\times(n-r)$ submatrices.
Let $\cF$ be the set of those $r$-sets $F\subset [n]$ where the 
 columns of $A$ labeled by the elements of $[n]\setminus F$ have full rank. 
We have  $|\cF|> c_0 \binom{n}{r}$. 
 
Let $\mathcal S$ be the $(n-r)$-dimensional subspace generated by the rows of $A$ 
  in ${\mathbf F}_2^n$ and let ${\mathcal R}$ be a subspace of dimension $r$ 
  such that ${\mathcal S}+{\mathcal R}$ is the whole space. 
Decompose the $n$-dimensional space into $2^r$ disjoint affine subspaces
 $$
   {\mathbf F}_2^n = \bigcup_{{\mathbf v}\in {\mathcal R}} ({\mathcal S}+{\mathbf v}). 
   $$
For any set $F\subset [n]$ let $\widehat{F}$ be a $0$-$1$ vector with support $F$. 
For each ${\mathbf v}\in {\mathcal R}$ let
$$
   {{\mathcal F}({\mathbf v})}:= \{ F: F\in \cF, \, \widehat {F}\in ({\mathcal S}+{\mathbf v})\}.
   $$
We have partitioned $\cF$ into  $2^r$ pairwise disjoint $r$-uniform families. 
Given any $F\in \cF$ the vectors of $({\mathcal S}+{\mathbf v})$ truncated to $([n]\setminus F)$ 
 are all distinct. 
Hence each ${\mathcal F}({\mathbf v})$ is a cancellative family. \qed

There are $\Theta(r^2)$ non isomorphic hypergraphs consisting of three edges
 $\{A, B, C\}$ with $A\cup B=A\cup C$. 
The {\it Tur\'an number} of the class of $r$-uniform hypergraphs 
 $\bbH:=\{ \cH_1, \cH_2, \dots \}$  is denoted by $\ex(n,\bbH)$. 
It is the size of the largest $r$-graph on $n$ vertices avoiding every $\cH\in \bbH$.
The sequence $\ex(n,\bbH )\binom{n}{r}^{-1}$ is monotone decreasing, 
 its limit is denoted by $\pi(\bbH)$. 
When we consider the determination of $c_r(n)$ as a Tur\'an type problem, 
 then there is a score of forbidden hypergraphs.
Take only one of them, namely $\bbG_r^3$ defined by three sets on $2r-1$ elements
 $[r]:=\{1,2,\dots, r\}$, $[r-1]\cup \{r+1 \}$ and $\{ r, r+1, \dots, 2r-1\}$.  
 It was proved in~\cite{Fur76} that 
 $$\left( \binom{r}{2} e^{1+ 1/(r-1)}\right)^{-1}\leq \pi(\bbG_r^3)\leq 
  \left( e\binom{r-1}{2}\right)^{-1}. 
    $$
 Concerning another case, for an even $r$ when $\bbT_r$
  is a blown up triangle, its three edges are $X\cup Y$, $Y\cup Z$, and $Z\cup X$
  where $|X|=|Y|=|Z|=r/2$ Frankl~\cite{Frankl} and Sidorenko~\cite{Sid_3, Sid_92}
 showed independently that $\pi (\bbT_r)=1/2$. 
More on this see~\cite{KeeSud05}.

\section{Conclusion, problems}

One of our  main results is 
 to give a better upper bound for the size of $2$-cancellative codes.
We {\bf conjecture} that the upper bounds of Theorems~\ref{th:2c} and~\ref{th:t*canc}
 are much closer to the truth than the simple probabilistic lower bounds we have.
This is probably true for the uniform case, see (\ref{eq:up_odd}), too.  

\begin{conjecture}\quad 
 $ n^{k+1-o(1)}< c_{2k+1}(n,2)= o(n^{k+1})$ as $n\to \infty$ and $k$ is fixed.
  \end{conjecture}

Call a code $\cF$ $t^*$-cancellative if 
 $$A_1\cup\dots \cup A_t\cup  B = A_1\cup \dots A_t \cup C \enskip  \Longrightarrow\enskip  
  B=C \,\text{ or } \,  \{ B, C\} \subset \{ A_i, \dots, A_t\}$$
for every  $t+2$ member {\it sequence} from $\cF$, and let
 $c^*(n,t)$ be the maximum size of such a code $\cF\subset 2^{[n]}$.
Obviously $C(n,t)\leq c^*(n,t)\leq C(n,t+1)\leq c(n,t)$.
One wonder if equality holds in some of these, and what other 
 relations these functions can have.

Using the Erd\H os, Frankl, and R\"odl~\cite{EFR} estimate, see (\ref{eq:EFR}), one has
 $$
   n^2 e^{-\alpha_r \sqrt{\log n}} \leq f_r(n,3(r-2)+3,3) \leq c_r(n,r-1).
   $$
The general upper bound (\ref{eq:gen_rt}) for $c_r(n,t)$ here only gives $O(n^3)$, but 
 in this case, leaving out those $r$-sets having an own pair, one can easily prove
$$
 c_r(n,r-1)\leq \binom{n}{2}. 
   $$
More of these type of problems, see, e.g.,~\cite{FR}.  

In Section~\ref{ss:proof_of_2k} the
 $2k$-partite hypergraph $\cF$ (with partite sets $V_1, \dots, V_{2k}$) 
  has an interesting property.
For every three members $A,B,C$ there exits a class $V_i$
 such that $A\cap V_i$, $B\cap V_i$ and $C\cap V_i$ are distinct.
It is natural to ask what other small substructures can be avoided this way.

The proof of Theorem~\ref{th:tcanc} concerning $c(n,t)$ presented in 
 Section~\ref{ss:tcanc} actually gives a slightly better upper bound.
With a little more calculations one can obtain an explicit bound $\kappa_t$ for $t\geq 3$
 such that 
 $$
  \limsup_n \left( c(n,t)\right)^{1/n} \leq \kappa_t < \frac{t+3}{t+2}. 
   $$
   
Many problems remain open.

\small

\bibliographystyle{amsalpha}

\end{document}